\documentstyle[11pt]{article}
\font\tenmsb=msbm10
\font\sevenmsb=msbm7
\font\fivemsb=msbm5
\newfam\msbfam
 \textfont\msbfam=\tenmsb
 \scriptfont\msbfam=\sevenmsb
 \scriptscriptfont\msbfam=\fivemsb
\def\Bbb#1{{\fam\msbfam #1}}
\def\L{{\cal L}}
\def\O{{\cal O}}
\def\P{{\cal P}}
\def\Q{{\cal Q}}
\def\N{{\cal N}}
\def\I{{\cal I}}

\def\E{{\cal E}}

\def\bZ{{\bf Z}}
\def\p{{\Bbb P}}
\def\Aone{{\Bbb A}^1}
\def\Pthree{{\Bbb P}^3}

\def\half{{1 \over 2}}

\def\mod{\mathop{\rm mod}}
\def\depth{\mathop{\rm depth}}
\def\supp{\mathop{\rm Supp}}
\def\rank{\mathop{\rm rank}}
\def\cok{\mathop{\rm coker}}
\def\ker{\mathop{\rm ker}}
\def\and{\mathop{\rm  and }}

\def\ra{\rightarrow}

\newtheorem{ther}{Theorem}[section]

\newtheorem{prop}[ther]{Proposition}
\newtheorem{cor}[ther]{Corollary}
\newtheorem{defn}[ther]{Definition}
\newtheorem{rmk}[ther]{Remark}
\newtheorem{ex}[ther]{Example}
\newenvironment{pf}
   {\noindent {\bf Proof:}}{\vskip .20in}

\title{The Hilbert Schemes of Degree Three Curves are Connected}
\author{Scott Nollet \\
        University of California, Berkeley \\
        Berkeley, CA 94720 \\
        {\normalsize{\em e-mail}: nollet@math.berkeley.edu}}
\date{}
\begin{document}
\maketitle
\begin{abstract} In this paper we show that the Hilbert scheme $H(3,g)$ of
locally Cohen-Macaulay curves in $\Pthree$ of degree three and genus
$g$ is connected. This is achieved by giving a classification of
these curves, determining the irreducible components of $H(3,g)$, and
giving certain specializations to show connectedness. As a byproduct, we
find that there are curves which lie in the closure of each irreducible
component.
\end{abstract}

\setcounter{section}{-1}

\section{Introduction}

In his thesis \cite{HC}, Hartshorne showed that the (full) Hilbert scheme
for projective subschemes with a fixed Hilbert polynomial is connected.
Often one studies certain subsets of the Hilbert scheme which parametrize
subschemes satisfying a certain property. For example, one can consider the
Hilbert scheme of smooth curves in $\Pthree$. The smooth curves of degree $9$
and genus $10$ afford an example in which this Hilbert scheme is {\it not}
connected (see \cite{H}, IV, ex. 6.4.3). It is not known for which properties
the corresponding Hilbert scheme is connected.
\par
In the present paper, we are interested in the Hilbert scheme $H(d,g)$
of locally Cohen-Macaulay curves in $\Pthree$ of degree $d$ and arithmetic
genus $g$. By work of several authors \cite{GSP,OS,S}, it is known that
$H(d,g)$ is nonempty precisely when $d>0, g=\half(d-1)(d-2)$ or
$d > 1, g \leq \half(d-2)(d-3)$. In a recent paper \cite{MDP3},
Martin-Deschamps and Perrin prove that $H(d,g)$ is
reduced only when $d=2$ or $g=\half(d-1)(d-2)$ or $g=\half(d-2)(d-3)$ or
$(d,g)=(3,-1)$. For all other $(d,g)$ pairs for which $H(d,g)$ is nonempty,
there is a nonreduced irreducible component corresponding to curves which
are extremal in the sense that their Rao modules have the largest possible
dimension.
\par
It can be gleaned from several sources \cite{MDP3,E1,E2} that $H(d,g)$ is
irreducible precisely in the cases $d=2$ or $g > \half(d-3)(d-4)+1$ or
$(d,g) \in \{(4,1),(3,1),(3,0),(3,-1)\}$. In particular, $H(d,g)$ is
connected in these cases. In the present paper we show that
$H(3,g)$ is connected for all $g$. This is the first interesting
case in the sense that these Hilbert schemes have several irreducible
components. Curiously, there certain extremal curves lie in the closure
of each irreducible component. \par
The paper is organized as follows.  In the first section, we review
several results of Banica and Forster \cite{BF} on multiplicity structures
on smooth curves in a smooth threefold and classify space curves
of degree two as an example. We also briefly review the extremal curves
studied by Martin-Deschamps and Perrin. In the second section, we
classify the multiplicity three structures on a line. This is used in
the third section, where we classify all locally Cohen-Macaulay curves of
degree $3$ in $\Pthree$. In particular, the irreducible components
of the Hilbert scheme are determined. We also produce some flat families of
triple lines, which show that the Hilbert scheme is connected.
\par
In this paper, we work over an algebraically closed field $k$ of arbitrary
characteristic. $S=k[x,y,z,w]$ denotes the homogeneous coordinate ring of
$\Pthree_{k}$. If $V \subset S$ is a closed subvariety, then $S_V$ denotes
the the homogeneous coordinate ring $S/I_V$ of $V$. We often use the
abbreviation CM to mean locally Cohen-Macaulay. $H(d,g)$ denotes the
Hilbert scheme of locally Cohen-Macaulay curves in $\Pthree$ of degree
$d$ and arithmetic genus $g$.

\section{Preliminaries}

In this section we review the results of Banica and Forster \cite{BF} on
multiplicity structures on smooth curves in smooth threefolds. As an
example, we give the classification of double lines in $\Pthree$, which
will be used in section two when we classify the triple lines in $\Pthree$.
We also recall a few notions from linkage theory and summarize the results
of Martin-Deschamps and Perrin on extremal curves.
\begin{defn}{\em
If $Y$ is a scheme, then a {\it locally Cohen-Macaulay multiplicity
structure} $Z$ {\it on} $Y$ is a locally Cohen-Macaulay scheme $Z$
which contains $Y$ and has the same support as $Y$. For short, we
simply say that $Z$ is a {\it multiplicity structure} on $Y$.
             \em}
\end{defn}
In \cite{BF}, Banica and Forster consider a smooth curve $Y$ inside
a smooth threefold $X$. Starting with a multiplicity structure
$Y \subset Z \subset X$, they define a filtration on $Z$ as follows.
Let $Y^{(i)}$ denote the subscheme of $X$ defined by ${I_Y}^i$.
Let $Z_i$ denote the subscheme of $X$ obtained from $Z \cap Y^{(i)}$
by removing the embedded points. This gives the (unique) largest
Cohen-Macaulay subscheme contained in $Z \cap Y^{(i)}$. If $k$
is the smallest integer such that $Z \subset Y^{(k)}$, we obtain the
{\it Cohen-Macaulay filtration} for $Y \subset Z$
$$Y=Z_1 \subset Z_2 \subset \dots \subset Z_k=Z.$$
\par
Letting $\I_i=\I_{Z_i}$, there are sheaves $L_j = \I_j/{\I_{j+1}}$
associated to this filtration. For any $i,j \geq 1$, it turns out that
${\I_i}{\I_j} \subset \I_{i+j}$, and hence the $L_j$ are $\O_Y$-modules.
In fact, the $L_j$ are shown to be locally free $\O_Y$-modules. Further,
there are induced multiplication maps $ L_i \otimes L_j \ra L_{i+j} $, which
are generically surjective (because $\I_j = \I_Z + {\I_Y}^i$ on an open set).
In particular, we get generically surjective maps $L_1^{\otimes j} \ra L_j$.
\par
{}From the above, we see that if $Z$ is a multiplicity structure on $Y$,
then there is a filtration $\{Z_j\}$ and exact sequences
$$0 \ra \I_{Z_{j+1}} \ra \I_{Z_j} \ra L_j \ra 0$$
where the $L_j$ are vector bundles on $Y$. If $Y$ is connected, the
multiplicity of $Z$ is defined by $\mu(Z)=\dim_K(\O_{Z,\eta})$, where
$\eta$ is the generic point of $Y$ and $K=\O_{Y,\eta}$ is the function
field of $Y$. The sequences above show that
$\mu(Z)=1+\sum_{j=1}^{k}{\rank L_j}.$
\begin{rmk}{\em
The above constructions can be carried out in $Z$ (instead of $X$),
and would yield the same filtration as above. Thus the Cohen-Macaulay
filtration is well-defined for abstract (non-embedded) multiplicity structures.
            \em}
\end{rmk}
Now we use the fact that $X$ is a smooth threefold. In this case, the
conormal sheaf $\I_Y/{\I_Y^2}$ is a rank two bundle on $Y$. Since the
surjection $\I_Y \ra L_1$ factors through the conormal bundle, we see
that $L_1$ has rank zero, one, or two. If the rank is zero,
then $L_1=0$ and the generically surjective maps show that all the $L_j=0$,
hence $Z=Y$. If the rank is two, then the surjection $\I_Y/{\I_Y^2} \ra L_1$
becomes an isomorphism, hence $Y^{(2)} \subset Z_2$ and $Z$ has generic
embedding dimension three. \par
We are mainly interested in the case rank$(L_1)=1$, in which case we
say the extension $Y \subset Z$ is {\it quasi-primitive}. In this case, the
generically surjective maps $L^{\otimes j} \ra L_j$ show that there are
divisors $D_j$ on $Y$ such that $L_j=L^j(D_j)$. The multiplication maps
$L_i \otimes L_j \ra L_{i+j}$ show that $D_i + D_j \leq D_{i+j}$ for all
$i,j \geq 1$ (define $D_1=0$). We say that $(L,D_2,\dots D_k)$
is the {\it type} of the extension $Z$. \par
\begin{rmk}{\em
The condition that $Y \subset Z$ be quasi-primitive is equivalent to
the condition that the generic embedding dimension of $Z$ is $2$.
If $z \in Z$ is a point where the embedding dimension is $2$, there is
an open neighborhood $U$ about $z$ and a smooth surface $S \subset U$ such that
$Z \subset S$. $Y$ is Cartier on $S$. If $t$ is a local equation for $Y$,
then $t^i$ gives a local equation for $Z_i$ on $S$.
            \em}
\end{rmk}
Having reviewed the theory of multiplicity structures on curves, we
present as an example the simplest case, the double structures on a
line in $\Pthree$.
\begin{prop}\label{2line}
Let $Y$ be the line $\{x=y=0\}$ in $\Pthree$ and let $a \geq -1$
be an integer. Let $f$ and $g$ be two homogeneous polynomials of degree
$a+1$ which have no common zeros along $Y$. Then $f$ and $g$ define a
surjection $u: \I_Y \ra \O_Y(a)$ by $x \mapsto f, y \mapsto g$. The
kernel of $u$ gives the ideal sheaf of a multiplicity two structure $Z$
on $Y$. Further, we have \\
(a) $p_a(Z)=-1-a$ \\
(b) $H^1_*(\I_Z) \cong (S/(x,y,f,g))(a)$. \\
(c) $I_Z=(x^2,xy,y^2,xg-yf)$. \\
(d) If $f^\prime,g^\prime$ define another two structure $Z^\prime$, then
$Z=Z^\prime$ if and only if there exists $c \in k^*$ such that
$f^\prime=cf \mod I_Y$ and $g^\prime=cg \mod I_Y$. \\
(e) Each multiplicity two structure $Z$ on $Y$ arises by the construction
above. \\
\end{prop}
\begin{pf}
This can be found in work of Migliore \cite{M} and also by work of
Martin-Deschamps and Perrin (\cite{MDP1}, IV, example 6.9) in the context
of linkage theory.
{}From the above theory of multiplicity structures, we see that giving a double
structure $Z$ on $Y$ is equivalent to finding a surjection $u:\I_Y \ra \L$,
where $\L$ is a line bundle on $Y$. Since such a surjection must factor
through $\I_Y/{\I_Y^2} \cong \O_Y(-1)^2$, we see that $\L \cong \O_Y(a)$
with $a \geq -1$ and that the map is given by two polynomials $f,g$ of
degree $a+1$.
\end{pf}
\begin{rmk}\label{smooth}{\em
If $Z$ is a double line from proposition \ref{2line} above, the exact sequence
$$0 \ra \O_Y(a) \ra \O_Z \ra \O_Y \ra 0$$
shows that $Z$ has local embedding dimension two at each point. In fact, $Z$ is
contained in a smooth (global) surface of degree $a+2$. To see this, one
can choose the polynomials $f,g$ in the variable $z,w$. Since $f$ and
$g$ have no common zeros along $Y$, the surface with equation $xg-yf$ is
smooth along $Y$. When $a=-1$,
this surface is a (smooth) plane which contains $Y$. When $a \geq 0$,
we note that the general surface of degree $a+2$ which contains $Y^{(2)}$ is
smooth away from $Y$. If we intersect these open conditions in
${\p}H^0(\I_Z(a+2))$, we find that there are surfaces of degree $a+2$
containing $Z$ which are smooth everywhere. The general surface containing
$Z$ of higher degree will have a finite number of singularities along $Y$.
            \em}
\end{rmk}
\begin{cor}\label{2linehilb}
Description of $H(2,g)$: \\
(a) If $g>0$, then $H(2,g)$ is empty. \\
(b) If $g=0$, then $H(2,g)$ is irreducible of dimension 8. All curves
in $H(2,g)$ are planar, and the general member is a smooth conic.
The reduced reducible curves (two lines meeting at a point) form
an irreducible family of dimension 7, and the multiplicity two structures
on a line form an irreducible family of dimension 5. \\
(c) If $g=-1$, then $H(2,g)$ is irreducible of dimension 8. The general curve
is a union of two skew lines. The multiplicity two
structures on a line form an irreducible family of dimension 7. \\
(d) If $g<-1$, then $H(2,g)$ is irreducible of dimension 5-2g. all
curves are multiplicity two structures on a line with $a=-1-g$ \\
\end{cor}
\begin{pf}
(a) is known, since ${\half}(d-1)(d-2)$ is an upper bound on the genus of
locally Cohen-Macaulay curves (see \cite{GSP} or \cite{OS}). The descriptions
of the families of reduced curves is standard. To describe the moduli for the
double lines of genus $g \leq 0$, we use proposition \ref{2line}. The choice
of the line $Y$ is given by a $4$-dimensional (irreducible) Grassman variety.
Given the line $Y$, the multiplicity structure $Z$ is uniquely determined
by the open set of $(f,g) \in H^0(\O_Y(a+1))^2/{k^*}$ where $f$ and $g$ have
no common multiple. This is an irreducible choice of dimension $2a+3 = 1-2g$.
Adding the choice of the line $Y$ gives an irreducible family of dimension
$5-2g$.
\end{pf}
In their excellent book \cite{MDP1}, Martin-Deschamps and Perrin build
a strong foundation for linkage theory of locally Cohen-Macaulay
curves in $\Pthree$.
Perhaps the most important result there is the structure theorem
for even linkage classes (see \cite{MDP1},IV,theorem 5.1). It states that
if $\L$ is an even linkage class of curves which are not arithmetically
Cohen-Macaulay, then there exists a curve $C_0 \in \L$ such that any other
curve $C \in \L$ is obtained from $C_0$ be a sequence of ascending elementary
double links (see \cite{MDP1}, III, definition 2.1) followed by a deformation
with constant cohomology through curves in $\L$. In particular, $C_0$ achieves
the smallest degree and genus among curves in $\L$. $C_0$ is called a
{\it minimal curve} for $\L$.
\par
A practical aspect of \cite{MDP1} (see chapter IV) is an algorithm for
finding a minimal curve associated to a finite length graded $S$-module.
If $M$ is such a module, there exists a minimal graded free resolution
$$0 \ra L_4 \stackrel{\sigma_4}{\ra} L_3 \stackrel{\sigma_3}{\ra}
L_2 \stackrel{\sigma_2}{\ra} L_1 \stackrel{\sigma_1}{\ra} L_0 \ra M \ra 0$$
which sheafifies to an exact sequence of direct sums of line bundles.
Let $\N_0=\ker \tilde{\sigma_2}$ and $\E_0=\ker \tilde{\sigma_3}$. The
algorithm of Martin-Deschamps and Perrin gives a way to split $\L_2$ into
$\P \oplus \Q$, where $\P$ and $\Q$ are also direct sums of lines bundles.
If $C_0$ is a minimal curve, then there exists an integer $h_0$ and exact
sequences
$$0 \ra \P \ra \N_0 \ra \I_{C_0}(h_0) \ra 0$$
$$0 \ra \E_0 \ra \Q \ra \I_{C_0}(h_0) \ra 0.$$
One consequence of this is that if
$$0 \ra L_4 \stackrel{\theta}{\ra} L_3 \ra Q \ra I_{C_0} \ra 0$$
is a minimal graded $S$-resolution for the total ideal of the minimal curve
$C_0$, then $L_3^\vee \stackrel{\theta^\vee}{\ra} L_4^\vee$ begins a minimal
resolution for $M^*$. \par
In a later paper \cite{MDP2}, Martin-Deschamps and Perrin
tackled the problem of bounding the dimensions of the Rao module of a curve
in terms of the degree and genus. For a curve $C$, define the
{\it Rao function} $\rho_C$ by $\rho_C(n)=h^1(\I_C(n))$. $r_a$ (resp. $r_o$)
is the smallest (resp. largest) value $n$ for which $\rho_C(n) \neq 0$.
Their bounds can be stated as follows (see \cite{MDP2}, theorem 2.5 and
corollary 2.6).
\begin{prop}\label{raobound}
Let $C \subset \Pthree$ be a curve of degree $d$ and genus $g$. Set
$$l=d-2,a=\half(d-2)(d-3)-g.$$
Then $a \geq 1,l \geq 0$, and the Rao function is bounded by \\
(1) $r_a \geq -a+1$. \\
(2) $\rho_C(n) \leq a$ for $0 \leq n \leq l$.\\
(3) $r_o \leq a+l-1$. \\
\end{prop}
The question of sharpness for these bounds has a nice answer. Not only can
equalities in (1), (2) and (3) be realized individually, but they can be
realized by one curve. Martin-Deschamps and Perrin call such a curve
{\it extremal}. These curves are characterized in theorem \ref{extremal}
below (see \cite{MDP3}, proposition 0.6 and theorem 1.1). In theorem
\ref{nilcomp} which follows (see \cite{MDP3}, theorem 4.2 and theorem 4.3),
it is shown that for $d \geq 3$, these curves form a nonreduced irreducible
component of $H(d,g)$. A finite length graded $S$-module $M$ is said to be a
{\it Koszul module parametrized by} $a \geq 1$ {\it and} $l \geq 0$ if
$M$ is isomorphic to a complete intersection module $S/(f_1,f_2,f_3,f_4)$
with $\deg f_1 = \deg f_2 = 1$, $\deg f_3=a$ and $\deg f_4=a+l$.
\begin{ther}\label{extremal}
Characterization of extremal curves:\\
(a) Fix $a \geq 1,l \geq 0$, and let $M$ be a Koszul module parametrized by
$a$ and $l$. Then any minimal curve for the even linkage class $\L(M)$ is
an extremal curve of degree $d=l+2$ and genus $g=-a+\half(d-2)(d-3)$.\\
(b) Conversely, let $C \subset \Pthree$ be an extremal curve of degree
$d \geq 2$ and genus $g < \half(d-2)(d-3)$. If $l=d-2$ and
$a=\half(d-2)(d-3)-g$, then $C$ is a minimal curve for a Koszul module
parametrized by $a$ and $l$.
\end{ther}
\begin{ther}\label{nilcomp}
Let $d \geq 3$ and $g < \half(d-2)(d-3)$. Then the family of extremal
curves gives an irreducible component of the Hilbert scheme
$H(d,g)$ of dimension ${3 \over 2}d(d-3)+9-2g$. This component is
nonreduced except when $(d,g)=(3,-1)$.
\end{ther}

\section{Triple Lines in $\Pthree$}

In this section we classify the multiplicity three structures on a fixed
line $Y \subset \Pthree$. If $W$ a quasiprimitive multiplicity three
structure of type $(L,D_2)$, then we have two exact sequences
\begin{equation}\label{filt1}
0 \ra \I_Z \ra \I_Y \ra \O_Y(a) \ra 0
\end{equation}
\begin{equation}\label{filt2}
0 \ra \I_W \ra \I_Z \ra \O_Y(2a+b) \ra 0
\end{equation}
where $Z$ is one of the multiplicity two structures on $Y$ described in
proposition \ref{2line}, $L = \O_Y(a), a \geq -1, \deg D_2=b \geq 0$.
We loosely say that $W$ is of {\it type} $(a,b)$.
\par
In classifying the triple lines of type $(a,b)$, we will handle the case
$a = -1$ separately. This is because the corresponding double line $Z$ is a
complete intersection when $a=-1$, while this is not the case for $a > 0$.

\begin{prop}\label{constr1}
Let $Y \subset \Pthree$ be the line
$\{x=y=0\}$ and let $Z$ be the multiplicity two structure $\{x=y^2=0\}$ on
$Y$. Let $p,q$ be two homogeneous polynomials of degrees $b-1,b$ which have
no common zeros along $Y$. Then $p$ and $q$ define a surjection
$u:\I_Z \ra \O_Y(b-2)$ by $x \mapsto p, y^2 \mapsto q$. The kernel of $u$
is the ideal sheaf of multiplicity three structure $W$ on $Y$.
Further, we have \\
(a) $p_a(W)=1-b$ \\
(b) $H^1_*(\I_W) \cong (S/(x,y,p,q))(b-2)$ \\
(c) $I_W=(x^2,xy,y^3,xq-y^2p)$ \\
(d) If $p^\prime,q^\prime$ define another three structure $W^\prime$, then
$W=W^\prime$ if and only if there exists $c \in k^*$ such that
$p^\prime=cp \mod I_Y$ and $q^\prime=cq \mod I_Y$. \\
(e) $W$ is quasiprimitive with second CM filtrant $Z$, unless
$b=1$ and $q=0$, in which case $W=Y^{(2)}$. \\
\end{prop}
\begin{pf}
Since $Z$ is a global complete intersection with total ideal $(x,y^2)$,
$I_Z \otimes S_Y = I_Z/{I_ZI_Y} \cong S_Y(-1) \oplus S_Y(-2)$
is freely generated by the images of $x$ and $y^2$. The map
$x \mapsto {\overline p}, y^2 \mapsto {\overline q}$ defines
a graded homomorphism $\phi:I_Z \ra I_Z/{I_ZI_Y} \ra S_Y(b-2)$. Since
$({\overline p},{\overline q})$ form a regular sequence in $S_Y$, the kernel
of the map $I_Z/{I_ZI_Y} \ra S_Y(b-2)$ is given by the Koszul relation
$x{\overline q}-y^2{\overline p}$, hence
$\ker \phi=(I_ZI_Y,xq-y^2p)=(x^2,xy,y^3,xq-y^2p)$. The cokernel
$\cok \phi=S_Y(b-2)/({\overline p},{\overline q}) \cong (S/(x,y,p,q))(b-2)$
has finite length, hence $\phi$ sheafifies to a surjection
$u: \I_Z \ra \O_Y(b-2)$. \par
Letting $W$ be the subscheme defined by $\I_W=\ker u$, we have an exact
sequence
$$0 \ra \I_W \ra \I_Z \ra \O_Y(b-2) \ra 0.$$
Since $H^0_*(u)=\phi$ and $H^1_*(\I_Z)=0$, we immediately deduce properties
(b) and (c). The snake lemma provides a second exact sequence
$$0 \ra \O_Y(b-2) \ra \O_W \ra \O_Z \ra 0,$$
which shows that $W$ is supported on $Y$ and that $\depth \O_W \geq 1$,
hence $W$ is a CM multiplicity three structure on $Y$ with $p_a(W)=1-b$.
If the polynomials $p^\prime,q^\prime$ also define $W$ by the construction
above, then $(q,p)$ and $(q^\prime,p^\prime)$ generate the same principal
$S_Y$-submodule of $S_Y(-1) \oplus S_Y(-2)$, hence property (d) holds. \par
If $q=0$, then $p$ must be a unit, as otherwise $p$ and $q$ will have common
zeros along $Y$. In this case, we see that $b-1 = \deg p =0$ and that
$I_W=I_Y^2$, whence $W=Y^{(2)}$. If $q \neq 0$, we use part (c) to see that
$I_W+I_Y^2=(x^2,xy,y^2,xq)=(I_Y^2,xq)$. At the points on
$Y$ where $q \neq 0$, this ideal is simply $(y^2,x)$, so the cokernel
of the inclusion $(I_Y^2,xq) \subset (y^2,x)$ has finite support. Since
the latter ideal defines the multiplicity two structure $Z$ on $Y$, we see
that $Z$ is the second CM filtrant for $W$, proving part (e).
\end{pf}
\begin{cor}\label{3line-1}
Triple lines of type $a=-1,b \geq 0$: Let $W$ be a quasiprimitive
multiplicity three structure on a line $Y$ of type $(-1,b)$ or the second
infinitesimal neighborhood $Y^{(2)}$. Then, after a suitable change of
coordinates, $W$ is constructed by proposition \ref{constr1}. The family of
such multiplicity three structures is irreducible of dimension $5+2b$.
\end{cor}
\begin{pf}
If $W=Y^{(2)}$ is given by the construction in taking $b=1,q=0,p=1$.
If $W$ is quasiprimitive, then we have the exact sequence \ref{filt2},
and the construction above gives all surjections $u:\I_Z \ra \O_Y(b-2)$.
To parametrize this family, we first choose the double line $Z$ (an
irreducible choice of dimension $5$ by corollary \ref{2linehilb}), and then
we must choose $(p,q) \in H^0(\O_Y(b-1)) \times H^0(\O_Y(b))/k^*$, which
is an irreducible choice of dimension $2b$. Thus the family is irreducible
of dimension $5+2b$.
\end{pf}
\begin{prop}\label{constr2}
Let $Z \subset \Pthree$ be the double line with total ideal
$I_Z=(I_Y^2,xg-yf)$, where $Y$ is the line $\{x=y=0\}$ and $f,g$ are
homogeneous polynomials of degree $a+1$ having no common zeros along $Y$,
as in proposition \ref{2line}. Let $p$ and $q$ be homogeneous polynomials of
degrees $b,3a+b+2$ having no common
zeros along $Y$. Then $p$ and $q$ define a surjection $u:\I_Z \ra \O_Y(2a+b)$
by $x^2 \mapsto pf^2, xy \mapsto pfg, y^2 \mapsto pg^2$ and
$xg-yf \mapsto q$. The kernel of $u$ is the ideal sheaf of a quasiprimitive
multiplicity three structure on $Y$ with second Cohen-Macaulay filtrant $Z$.
Further, we have \\
(a) $p_a(W)=-2-3a-b$ \\
(b) $I_W=(I_Y^3,x(xg-yf),y(xg-yf),p(xg-yf)-ax^2-bxy-cy^2)$, where $a,b,c$
are chosen so that $q=af^2+bfg+cg^2 \mod I_Y$.\\
(c) If $p^\prime,q^\prime$ define the multiplicity three structure
$W^\prime$, then $W=W^\prime$ if and only if there exists $d \in k^*$ such
that $p^\prime=dp \mod I_Y$ and $q^\prime=dq \mod I_Y$. \\
\end{prop}
\begin{pf}
The ideal $I_Z=(x^2,xy,y^2,xg-yf)$ has $S$-presentation
$$S(-a-3)^2 \oplus S(-3)^2 \stackrel{\varphi}{\ra} S(-a-2) \oplus S(-2)^3 \ra
I_Z \ra 0$$
given by the matrix
$$\varphi=\left(\begin{array}{cccc}
y & 0 & -g & 0 \\
-x & y & f & -g \\
0 & -x & 0 & f \\
0 & 0 & x & y \\
\end{array}\right).$$
Tensoring with $S_Y$, we see that $I_Z/{I_Z}{I_Y} = \cok {\varphi \otimes S_Y}$
is isomorphic to $S_Y(-a-2) \oplus (f^2,fg,g^2)(2a)$,
where ${\overline x}^2,{\overline x}{\overline y},{\overline y}^2$ are
identified with $f^2,fg,g^2$. Making this identification, we have an inclusion
$I_Z/{I_Z}{I_Y} \subset S_Y(-a-2) \oplus S_Y(2a)$ whose cokernel has finite
length. It follows that the sheafification of $I_Z/{I_Z}{I_Y}$ is isomorphic
to $\O_Y(-a-2) \oplus \O_Y(2a)$, freely generated by $xg-yf$ and an element
$e$ such that $ef^2={\overline x}^2,efg={\overline x}{\overline y}$ and
$eg^2={\overline y}^2$. \par
The polynomials $p$ and $q$ give a graded homomorphism
$$\phi:I_Z \ra I_Z/{I_Z}{I_Y} \subset S_Y(-a-2) \oplus S_Y(2a)
\stackrel{({\overline q},{\overline p})}{\longrightarrow} S_Y(2a+b).$$
The kernel of the map $({\overline q},{\overline p})$ is given by the
Koszul relation $qe-p(xg-yf)$. Since $f$ and $g$ are relatively prime of
degree $a+1$, the map $S_Y(-2a-2)^3 \ra S_Y$ given by $(f^2,fg,g^2)$ is
surjective in degrees $\geq 3(a+1)-1$, and hence there exist
$(a,b,c)$ such that $q=af^2+bfg+cg^2 \mod I_Y$ (because $\deg q \geq 3a+2$).
We can now write $(ax^2+bxy+cy^2-p(xg-yf))=I_Z/{I_Z}{I_Y} \cap \ker ({\overline
q},{\overline p})$, and hence
$$\ker \phi=(x^3,x^2y,xy^2,y^3,x(xg-yf),y(xg-yf),ax^2+bxy+cy^2-p(xg-yf)).$$
The cokernel of $\phi$ is of finite length, so $\phi$ sheafifies to a
surjection $u:\I_Z \ra \O_Y(2a+b)$.
\par
Letting $W$ be the subscheme whose ideal sheaf is the kernel of $u$, we get
an exact sequence
$$0 \ra \O_Y(2a+b) \ra \O_W \ra \O_Z \ra 0$$
which shows that $\supp W=Y$ and $\depth \O_W \geq 1$, hence $W$ is a
multiplicity three structure on $Y$. Since $p_a(Z)=-1-a$, the exact sequence
shows that $p_a(W)=-2-3a-b$. The exact sequence
$$0 \ra \I_W \ra \I_Z \ra \O_Y(2a+b) \ra 0$$ shows that $I_W=\ker \phi$.
If $p^\prime$ and $q^\prime$ define $W^\prime$ by the above construction and
$W^\prime=W$, then $eq^\prime-(xg-yf)p^\prime$ generates the same
$S_Y$-submodule of
$I_Z/{I_Z}{I_Y} \subset S_Y(-a-2) \oplus S_Y(2a)$ as $eq-(xg-yf)p$. Since
$e$ and $xg-yf$ are free generators, it follows that there exists
$d \in k^*$ such that $p^\prime=dp \mod I_Y$ and $q^\prime=dq \mod I_Y$.
This proves (a),(b) and (c). \par
{}From part (b), we find that $I_W+I_Y^2=(I_Y^2,h(xg-yf))$. The cokernel of
the inclusion $(I_Y^2,h(xg-yf)) \subset I_Z$ is supported on the zeros of
$h$ along $Y$. Since $Z$ has no embedded points, the second CM filtrant of $W$
is $Z$ and the extension $Y \subset W$ is quasi-primitive. \par
\end{pf}
\begin{rmk}\label{coh}{\em
Note that if $u:\I_Z \ra \O_Y(2a+b)$ is the surjection above, then
$H^1_*(u)$ is the zero map. Indeed, $H^1_*(\I_Z)$ is generated in degree
$-a$ by proposition \ref{2line} while $H^1(\O_Y(a+b))=0$. In particular,
we have an exact sequence
$$0 \ra \cok \phi \ra M_W \ra M_Z \ra 0$$
which shows that the Rao module $M_W$ is 2-generated. Since all curves of
degree $\leq 2$ have Rao modules which are zero or principal, $W$ is a
minimal curve. Further, the exact sequences
\ref{filt1} and \ref{filt2} now give that
$$h^2(\I_W(l))=h^1(\O_Y(l))+h^1(\O_Y(a+l))+h^1(\O_Y(2a+b+l))$$
for all $l \in \bZ$. \par
{}From the total ideal $I_W$ given in part (b), one can compute a minimal
graded $S$-resolution for $I_W$, which has the form
\begin{equation}\label{wres}
S(-a-b-4) \oplus S(-a-5)^2 \stackrel{\theta}{\ra} S(-a-b-3)^2 \oplus S(-a-4)^4
\ra S(-a-b-2) \oplus S(-3)^4
\end{equation}
This resolution determines $h^0(\I_W(l))$ for all $l \in \bZ$. Combining
with the dimensions $h^2(\I_W(l))$ found above, all the $h^i(\I_W(l))$ can
be computed. \par
The machinery for minimal curves of Martin-Deschamps and Perrin shows that
$\theta^\vee$ begins a minimal resolution for $M_W^*$. Completing this
resolution and dualizing the last map gives a presentation for $M_W$.
Carrying this out (we suppress the calculation here), one finds that
$H^1_*(\I_W) \cong \cok \psi$, where
$$\psi: S(2a+b-1)^2 \oplus S(a+b-1)^2 \oplus S(a-1)^2 \ra S(2a+b) \oplus S(a)$$
is the map given by the matrix
$$\left(\begin{array}{cccccc}
x & y & fp & gp & -gc & -fa-gb \\
0 & 0 & x & y & f & g \\
\end{array}\right) $$
Here $a,b$ and $c$ are chosen as in part (b) of the proposition.
\em}\end{rmk}
\begin{rmk}{\em
In the case when $b=0$, ${\overline p}$ must be a unit. It follows that
the generators $x(xg-yf)$ and $y(xg-yf)$ are not needed for the total
ideal $I_W$ (see also \cite{B}, p. 24). In this case it is clear that $W$ is
the unique triple line supported on $Y$ and contained in the surface
defined by $ax^2+bxy+cy^2-p(xg-yf)$.
\em}\end{rmk}
\begin{cor}\label{3line0}
Triple lines of type $a,b \geq 0$: Let $W$ be a quasi-primitive multiplicity
three structure on a line $Y \subset \Pthree$ of type $(a,b)$ with
$a,b \geq 0$. Then $W$ arises from the construction of proposition
\ref{constr2} after a suitable change of coordinates. The family of these
triple lines is irreducible of dimension $10+5a+2b$.
\end{cor}
\begin{pf}
Since $W$ is of type $(a,b)$ with $a \geq 0$, there is an exact sequence
$$0 \ra \I_W \ra \I_Z \stackrel{u}{\ra} \O_Y(2a+b) \ra 0$$
where $Z$ is a double line of type $a \geq 0$. By Proposition \ref{2line},
we may change coordinates so that $I_Z=(x^2,xy,y^2,xf-yg)$ where $f,g$ are
homogeneous polynomials of degree $a+1$ with no common zeros along $Y$.
As in the proof of \ref{constr2} above, $\I_Z \otimes \O_Y \cong
\O_Y(2a) \oplus \O_Y(-a-2)$ freely generated by $e$ and the image of $xg-yf$,
where $ef^2=x^2,efg=xy$ and $eg^2=y^2$. From this it is clear that such a
map $u$ is given by homogeneous polynomials $p,q$ of degrees $b,3a+b+2$
which have no common zero along $Y$, and hence $W$ arises by the construction
of proposition \ref{constr2}.
\end{pf}

\section{The Hilbert Scheme}

In this section we describe the Hilbert scheme $H(3,g)$ of locally
Cohen-Macaulay curves of degree $3$ and arithmetic genus $g \leq 1$.
In particular, we classify all CM curves of degree $3$ and describe
the irreducible components of $H(3,g)$. We also show that certain
extremal curves lie in the closure of each irreducible component, hence
that $H(3,g)$ is connected. We begin with the curves of genus
$-1 \leq g \leq 1$, which have been described elsewhere. \par
\begin{prop}\label{irred}
For $-1 \leq g \leq 1$, the Hilbert scheme $H(3,g)$ is smooth and irreducible
of dimension $12$.
\end{prop}
\begin{pf}
$H(3,1)$ consists of arithmetically Cohen-Macaulay curves, hence
we can apply a theorem of Ellingsrud \cite{E2}. That $H(3,0)$ and $H(3,-1)$
are smooth and irreducible of dimension $12$ is part of \cite{MDP3}, theorem
4.1. $H(3,0)$ consists of arithmetically Cohen-Macaulay curves and
$H(3,-1)$ consists of extremal curves.
\end{pf}
For $g \leq -2$, the Hilbert scheme is not irreducible, and more work is
required to show connectedness. Our first task is to describe how the
unions of double lines and reduced lines fit in with the irreducible families
of triple lines.
\begin{prop}\label{famI}
Fix $g \leq -2$. Then \\
(a) The family of curves $W = Z \cup_{2P} L$ formed by taking the union
of a double line $Z$ with $p_a(Z)=g-1$ and a line $L$ which meets $Z$ in
a double point form an irreducible family of dimension $9-2g$. \\
(b) The family of curves $W = Z \cup_P L$ formed by taking the union
of a double line $Z$ with $p_a(Z)=g$ and a line $L$ which meets $Z$ in
a reduced point form an irreducible family of dimension $8-2g$. \\
(c) The family of curves $W$ which are triple lines of type $(-1,1-g)$
form an irreducible family of dimension $7-2g$. \\
Each curve above is an extremal curve, hence is a minimal curve for a
complete intersection module with parameters $l=1$ and $a=-g$. The families
(b) and (c) lie in the closure of the family (a).
\end{prop}
\begin{pf}
Let $W=Z \cup_{2P} L$ be a curve from family (a) above. After a change of
coordinates we may write $I_L=(x,z)$ and $I_Z=(x^2,xy,y^2,xg-yf)$.
We have an exact sequence
$$0 \ra \I_W \ra \I_Z \oplus \I_L \stackrel{\pi}{\ra} \I_{2P} \ra 0$$
where $2P=Z \cap L$ denotes the double point. Noting that
$I_L+I_Z=(x,z,y^2,yf)$ and that $I_{2P}=(x,z,y^2)$ ($2P$ is a complete
intersection), we see that $H_0^*(\pi)$ is surjective. Since $H^1_*(\I_{2P})$
vanishes in positive degree, we conclude that $r_o(W)=r_o(Z)=-g$ and that
$\rho_W(1)=-g$. $Z(x^2,y^2z)$ links $W$ to $W^\prime=Z^\prime \cup_{2Q} L$,
which is also from family (a). Applying the argument above and using the
isomorphism $M_{W^\prime}^* \cong M_W(1)$ shows that $r_a(W)=1+g$ and
$\rho_W(0)=-g$. Thus $W$ is extremal. \par
To parametrize this family of curves, one first chooses the double line
$Z$ (an irreducible choice of dimension $7-2g$ by corollary
\ref{2linehilb}, since $p_a(Z)=g-1$), then a point $P \in Z$ ($1$ parameter),
and finally a line $L$ through $P$ lying in the tangent plane to $Z$ at $P$
($1$ parameter). This shows that this family is irreducible of dimension
$9-2g$. \par
The argument for $W=Z \cup_P L$ is similar. We choose suitable
coordinates and write $I_L=(x,z), I_Z=(x^2,xy,y^2,xg-yf)$. We have an exact
sequence
$$0 \ra \I_W \ra \I_Z \oplus \I_L \stackrel{\pi}{\ra} \I_P \ra 0$$
where $P=Z \cap L$. Writing $I_{P}=(x,y,z)$ and $I_L+I_Z=(x,z,y^2,yf)$ we
see that $\dim \cok H^0_*(\pi(l))=1$ for $1 \leq l \leq -g=\deg yf-1$.
It follows again that $r_o(W)=r_o(Z)-1=-g$ and $\rho_W(1)=\rho_Z(1)-1=-g$.
$Z(x^2,y^2z)$ links $W$ to a curve $W^\prime$ from family (b), so we find
that $W$ is extremal. \par
To parametrize these curves, we first choose a double line $Z$ (an
irreducible choice of dimension $5-2g$, since $p_a(Z)=g$), and then
choose a general line $L$ which meets $Z$ ($3$ parameters). This shows
that (b) is an irreducible family of dimension $8-2g$. \par
If $W$ is a triple line of type $(-1,1-g)$, then from corollary \ref{3line-1}
we have $M_W=(S/(x,y,p,q))(-1-g)$ where $\deg p = 1-g$ and $\deg q = 2-g$.
It follows that $W$ is extremal. The family of such triple lines $W$ is
irreducible of dimension $7-2g$ by corollary \ref{3line-1}.
\par
By theorem \ref{nilcomp}, the scheme $H_{\gamma,\rho}$ of extremal curves
is irreducible of dimension $9-2g$ when $g \leq -2$. It follows that
family (a) gives the general member of the family, and that families (b)
and (c) lie in the closure.
\end{pf}
\begin{prop}\label{famII}
Fix $g \leq -2$. Then\\
(a) The family of curves $W = Z \cup L$ formed by taking the union of a
double line $Z$ with $p_a(Z)=g+1$ and a disjoint line $L$ form an irreducible
family of dimension $7-2g$. \\
(b) The family of curves $W$ which are triple lines of type $(0,-2-g)$ form
an irreducible family of dimension $6-2g$. \\
The curves above are all minimal, and each curve in family (b) is obtained
from curves in family (a) by a deformation which preserves cohomology.
\end{prop}
\begin{pf}
Let $W = Z \cup L$ be a curve from family (a) above. We begin by computing the
total ideal and cohomology for $W$. In suitable coordinates, we may write
$I_L=(z,w),I_Z=(x^2,xy,y^2,xg-yf)$ with $g,f \in k[z,w]$. In particular,
$xg-yf \in I_L$ and hence $J=((x,y)^2(z,w),xg-yf) \subset I_L \cap I_Z$.
One can compute that the minimal graded $S$-resolution of $J$ is of the form
$$S(g-3)\oplus S(-5)^2 \ra S(g-2)^2\oplus S(-4)^7 \ra S(g-1)\oplus S(-3)^6 $$
and hence $J$ is the total ideal for $W$. Comparing with resolution
\ref{wres} (with $a=0,b=-2-g$) shows that $W$ has the same Hilbert function
as a curve in family (b). Moreover, the exact sequence
$$0 \ra \I_W \ra \I_Z \oplus \I_L \ra \O \ra 0$$
shows that $h^2(\I_W(l))=h^2(\I_L(l))+h^1(\I_Z(l))=h^1(\O_L(l))+
h^1(\O_Y(l))+h^1(\O_Y(-g-2+l))$. This agrees with the second cohomology
dimensions found in remark \ref{coh} for triple lines of type $(0,-2-g)$,
hence the dimensions of the cohomology groups for families (a) and (b)
are the same. This same exact sequence also shows that $M_W$ is
2-generated. \par
The family (a) is parametrized by first choosing $Z$ (an irreducible
choice of dimension $3-2g$ by corollary \ref{2linehilb}) and then choosing
a general line $L$ ($4$ parameters), hence the family (a) is irreducible
of dimension $7-2g$. Family (b) is irreducible of dimension $6-2g$ by
corollary \ref{3line0}. These curves are minimal because they are of
degree three and their Rao modules are 2-generated. \par
Let $W$ be a triple line from family (b). By corollary \ref{3line0}, we
can change coordinates and write $I_W=((x,y)^3,xq,yq,hq-ax^2-bxy-cy^2)$, where
$q=xg-yf$ is a quadric surface containing the underlying second CM filtrant
$Z$. By remark \ref{smooth}, $q$ may be chosen to be the equation of a
smooth quadric $Q$. We may choose $z$ and $w$ so that $q=xz-yw$. \par
On the smooth quadric $Q$ the family of lines $L_t=Z(x+wt,y+zt)$ give a
flat family over $\Aone$ with $L_0=Y$. $D_t=L_t \cup Y$ forms a flat family
such that $D_0=Z$ is the double line $Z$ on $Q$ supported on $Y$, the second
Cohen-Macaulay filtrant of $W$. Writing this family as
$D \subset \Pthree \times \Aone \stackrel{\pi}{\ra} \Aone$, we see by
By Grauert's theorem, $\pi_*(\I_D(-g))$ is locally free on $\Aone$, hence
globally free. In particular, if $s_1 \in \I_{D_1}(-g)$ is
the equation of a smooth surface containing $D_1=L_1 \cup Y$, we can find a
section $s_t$ extending $s_1$ such that $s_0=hq-ax^2-bxy-cy^2$. \par
Now consider the family $C_t=S_t \cap (Y^{(2)} \cup L_t)$. Let $U \subset
\Aone$
be the open set where $C_t$ is locally Cohen-Macaulay. For $t \neq 0$, $C_t$
is the disjoint union of a double line on $Y$ and the line $L_t$. The
ideal of $C_t$ is given by $I_t=((x,y)^2(x+wt,y+zt),s_t)$. Note that
$x^2(y+zt)-xy(x+wt)=xqt \in I_t$ and similarly $yqt \in I_t$.
Flattening over $U$, we must add $xq$ and $yq$ to $I_t$. In particular, the
limit ideal $I_0$ contains $((x,y)^3,xq,yq,s_0)$, and hence gives $W$.
\end{pf}

\begin{prop}\label{g-2}
The Hilbert scheme $H(3,-2)$ consists of the following pair of irreducible
components:\\
(a) The irreducible family $H_{-1}$ of dimension $13$ from proposition
\ref{famI}.\\
(b) The closure $H_0$ of the irreducible family of sets of three disjoint
lines.
This closure is $12$-dimensional and contains the curves from proposition
\ref{famII}.
\end{prop}
\begin{pf}
Let $W \in H(3,-2)$ be a curve. Since $g<0$, $W$ is not integral, hence
is reducible or nonreduced. If $W$ is reduced, then $W$ is the union of
$3$ disjoint lines (because any union of a conic and line has $g \geq -1$),
hence lies in family (b). From corollary \ref{2linehilb}, any double line with
arithmetic genus $-1$ is a limit of disjoint unions of two lines. Adding
another disjoint line to this deformation, we see that the disjoint union
of a double line $Z$ of genus $-1$ and a reduced line $L$ lies in the
closure of the family of three disjoint lines. Thus family (b) is
irreducible. \par
If $W$ is not reduced, then $\deg \supp W <3$. If $\deg \supp W=1$, then
$W$ is a triple structure on a line $Y$ of arithmetic genus $-2$, which can
only have type $(-1,3)$ or $(0,0)$ ($W$ must be quasiprimitive, since
otherwise $W=Y^{(2)}$, which satisfies $p_a(W)=0$). By propositions \ref{famI}
and \ref{famII}, these lie in families (a) and (b) respectively. If
$\deg \supp W=2$, then $W$ is a union of a double line $Z$ and a reduced
line $L$ (the support of $W$ cannot be an irreducible conic, since a multiple
conic has degree $\geq 4$). $Z$ can meet $L$ in a scheme of length $0$, $1$,
or $2$, hence $W$ lies in family (a) or family (b). \par
Family (a) cannot lie in the closure of family (b) by reason of dimension.
Family (b) cannot lie in the closure of family (a) by reason of semicontinuity;
the curves in family (a) are extremal, while the curves in family (b) are not
(If $C$ is in family (b), one checks that $h^1(\I_C(-1))=0$).
\end{pf}
\begin{prop}\label{g-3}
Let $g \leq -3$. Then the Hilbert scheme $H(3,g)$ consists of the following
irreducible components: \\
(a) The irreducible family of dimension $9-2g$ from proposition \ref{famI},
which we now denote $H_{-1}$.\\
(b) The closure of the irreducible family of dimension $7-2g$ from
proposition \ref{famII}, which we now denote $H_0$.\\
(c) for each $0 < a < (-2-g)/3$, the closure of the irreducible family
$H_a$ of dimension $14-2g-a$ consisting of triple lines of type
$(a,-2-3a-g)$.\\
\end{prop}
\begin{pf}
Let $C \in H(3,g)$. Then $C$ is not integral because $g \leq -3$. If $C$ were
reduced, it would be a union of $3$ lines (these have genus $\geq -2$, hence
are ruled out) or the union of a conic and a line (which has $\geq -1$, hence
is ruled out). Thus $C$ is not reduced and $\dim \supp C < 3$.
If $C$ has support of degree $2$, the support cannot be irreducible, since
a multiplicity structure on a conic has degree at least $4$. Hence the
support of $C$ consists of two lines, and all possible configurations are
covered in families (a) and (b) above. If $C$ has support of degree $1$,
then $C$ is a triple line and corollaries \ref{3line-1} and \ref{3line0}
show that $C$ is among the families listed above. \par
Now we show that the $H_i$ are irreducible components. Let $-1 \leq i < j \leq
(-2-g)/3$. $H_i$ is not contained in the closure of $H_j$ because
$\dim H_i > \dim H_j$. On the other hand, semicontinuity shows that
$H_j$ is not contained in the closure of $H_i$. Indeed, from corollaries
\ref{3line-1} and \ref{3line0}, we see that the Rao
module for a triple line of type $(a,b)$ as a generator of minimal degree
$-2a-b$, and hence a minimal degree generator for the Rao module
of a curve in $H_a$ occurs in degree $g+2+a$. This shows that for
$C \in H_i$ we have $h^0(\O_C(g+2+i)) \neq 0$ while for $C \in H_j$ we have
$h^0(\O_C(g+2+i))=0$. Hence there can be no specialization from a family
of curves in $H_i$ to a curve in $H_j$.
\end{pf}

\begin{prop}\label{spec}
For each $a \geq 0$ and $b \geq 0$, there exists a flat family
$W \subset \Pthree \times \Aone$ whose general member $W_t$ is a triple
line of type $(a,b)$ for $t \neq 0$ and whose special member $W_0$ is
a triple line of type $(-1,3a+b+3)$.
\end{prop}
\begin{pf}
Consider the family defined by the ideal $I_t$ with generators
$$x^3,x^2y,xy^2,y^3,x(xz^{a+1}-tyw^{a+1}),y(xz^{a+1}-tyw^{a+1})$$
$$z^bt^2(xz^{a+1}-tyw^{a+1})-x^2w^{a+b}.$$
We flatten this family over $t$ by adding to the ideal those elements
which are multiples of $t$. Let $A,B,C$ denote the last three generators
given for the ideal. Then we must add
$$D=(w^{a+b}A+z^{a+1}C)/t=-xyw^{2a+b+1}+z^{a+b+1}t(xz^{a+1}-tyw^{a+1})$$
to the ideal. We must also add
$$E=(w^{2a+b+1}B+z^{a+1}D)/t=-y^2w^{3a+b+2}+z^{2a+b+2}(xz^{a+1}-tyw^{a+1})$$
to the ideal. Setting $t=0$, we find that the limit ideal $I_0$ contains
the generators
$$x^3,x^2y,xy^2,y^3,xyz^{a+1},x^2w{a+b},xyw^{2a+b+1},xz^{3a+b+3}-y^2w^{3a+b+2}.$$
It follows that the saturation of $I_0$ contains the ideal
$$(x^2,xy,y^3,xz^{3a+b+3}-y^2w^{3a+b+2}),$$
but this is the total ideal of a triple structure of type $(-1,3a+b+3)$ by
corollary \ref{3line-1}. On the other hand, corollary \ref{3line0}
shows that the ideal $I_t$ for $t \neq 0$ is the total ideal of a triple
line of type $(a,b)$. This gives the flat family $W$.
\end{pf}
\begin{rmk}{\em
The commutative algebra in the proof above was inspired by a geometric
example of Robin Hartshorne. He gave an example of a deformation of
three disjoint lines to a triple line of type $(-1,3)$ by deforming
the unique quadric containing the three lines to a double plane while
at the same time bringing the lines together.
\em}\end{rmk}
\begin{ther}
The Hilbert scheme $H(3,g)$ is connected is connected whenever it is nonempty.
\end{ther}
\begin{pf}
By proposition \ref{irred}, it suffices to consider the case $g \leq -2$.
In this case $H(3,g)$ has irreducible components $\{H_a\}_{a \geq -1}$ by
propositions \ref{g-2} and \ref{g-3}. Let $H_a$ be one of these components
with $a \geq 0$. Choosing $b=-2-3a-g$, proposition \ref{spec} gives a
family of triple lines whose general member lies in $H_a$ and whose special
member lies in $H_{-1}$.
\end{pf}
\begin{rmk}{\em
The proof of proposition \ref{spec} shows that a triple line $W$ with total
ideal $(x^2,xy,y^3,xz^{1-g}-y^2w^{-g})$ lies in the closure of each irreducible
component of $H(3,g)$.
\em}\end{rmk}
\begin{ex}{\em
Hartshorne has shown that the Hilbert scheme $H(4,0)$ is also connected.
Here we give an independent proof using the methods of this paper. $H(4,0)$
has two irreducible components (\cite{MDP3},$\S 4$): $H_1=$ the extremal
curves (these have Rao module of Koszul type parametrized by
$a=1$ and $l=2$) and $H_2=$ the curves with Rao module $k$ in degree $1$.
We will give a specialization from quadruple lines in $H_2$ to quadruple
lines in $H_1$.
\par
Let $Y$ be the line $\{x=y=0\}$ and $W$ be the planar triple line with total
ideal $I_W=(x,y^3)$. As in proposition \ref{3line-1}, a pair $(h,k)$ of
homogeneous polynomials of degrees $1$ and $3$ with no common zeros along
$Y$ determines a map $I_W \stackrel{\phi}{\ra} S_Y$ by $x \mapsto h,
y \mapsto k$ which sheafifies to a surjection $u$. $\ker u=\I_T$ defines a
multiplicity four line $Z_1$ such that $p_a(Z_1)=0$,
$I_{Z_1}=(x^2,xy,y^4,xk-y^2h)$, and $H^1_*(\I_{Z_1}) \cong S/(x,y,h,k)$.
It follows that $Z_1 \in H_1$.
\par
Letting $V$ be a quasiprimitive multiplicity three structure of type $(-1,1)$
on $Y$, proposition \ref{3line-1} shows that we may write
$I_V=(x^2,xy,xq-y^2)$, with $q \not\in I_Y$ ($p$ is unit in this case).
As in the proof of proposition \ref{3line0}, a pair $(f,g)$ of forms with
no common zeros along $Y$ determines a map $I_V \ra S_Y(-1)$ by
$x^2 \mapsto 0,xy \mapsto f,xq-y^2 \mapsto g$ which sheafifies to a
surjection $w$. $\ker w=\I_{Z_2}$ defines a multiplicity four line $Z_2$ such
that $Z_2$ such that $p_a(Z_2)=0$,$I_{Z_2}=(x^2,xy^2,xyq-y^3,gxy-f(xq-y^2))$
and $H^1_*(\I_{Z_2}) \cong S/(x,y,f,g))(-1)$, hence $Z_2 \in H_2$.
\par
Now consider the ideal
$$I_t=(x^2,xy^2,ty^3-xyz,xyw-tz(y^2t-xz))$$
in the ring $k[t][x,y,z,w]$. For $t \neq 0$, this gives the total ideal
of a curve in $H_2$ (see $Z_2$ above). Flattening over $t$, we add to this
ideal the multiples of $t$. Letting $A,B,C$ denote the last three generators
listed, we add
$$D=(wB+zC)/t=y^3w+xz^3-ty^2z$$
$$E=(zA+yB)/t=y^4$$
to $I_t$. Setting $t=0$, it follows that
$$(x^2,xy^2,xyz,xyw,y^3w+xz^3,y^4) \subset I_0$$
and hence $(x^2,xy,y^4,y^3w+xz^3) \subset {\overline I_0}.$ This ideal gives
a multiplicity four line in $H_1$ (see $Z_1$ above). This shows that
$H(4,0)$ is connected.
\em}\end{ex}
\begin{rmk}{\em
The results in this paper raise several questions: \\
(1) Can each locally Cohen-Macaulay curve $C \subset \Pthree$ be deformed
with constant cohomology to a quasiprimitive multiple line?\\
(2) Can each multiplicity structure on a line be deformed to an extremal
multiplicity structure on the same line? \\
(3) Is $H(d,g)$ connected for all $(d,g)$?\\
The answers are yes when $d=2$ and $d=3$. Positive answers to (1) and (2)
would give a positive answer to (3).
\em}\end{rmk}

\end{document}